\documentclass[12pt]{article}
\title{Two-step flag manifolds and the Horn conjecture}
\author{Kevin Purbhoo\\ University of British Columbia}

\usepackage[mathscr]{eucal}
\usepackage{latexsym}
\usepackage{amssymb}
\usepackage{amsthm}
\usepackage{amsmath}
\usepackage{epsfig}
\usepackage{graphics}

\newcommand{\frn}{\mathfrak{n}}
\newcommand{\frb}{\mathfrak{b}}
\newcommand{\frg}{\mathfrak{g}}

\newcommand{\frp}{\mathfrak{p}}
\newcommand{\CC}{\mathbb{C}}

\newcommand{\codim}{\mathop{\mathrm{codim}}}
\newcommand{\Hom}{\mathrm{Hom}}
\newcommand{\Image}{\mathop{\mathrm{Image}}}
\newcommand{\Stab}{\mathop{\mathrm{Stab}}}
\newcommand{\Fstand}{F^{\mathrm{std}}}
\newcommand{\revddots}{\reflectbox{$\ddots$}}
\newcommand{\zeroset}{\{0\}}
\newcommand{\zero}{{\rm `$0$'} }
\newcommand{\one}{{\rm `$1$'} }

\newcommand{\zeros}{{\rm `$0$'}s }
\newcommand{\ones}{{\rm `$1$'}s }
\newcommand{\twos}{{\rm `$2$'}s }
\newcommand{\modmod}{\big /}

\newcommand{\sqempty}{\multicolumn{1}{|c}{}}
\newcommand{\squp}{\multicolumn{1}{|c}{$^*$}}
\newcommand{\sqdown}{\multicolumn{1}{|c}{$_+$}}
\newcommand{\squpdown}{\multicolumn{1}{|c}{$^*_+$}}

\newcommand{\squpend}{\multicolumn{1}{|c|}{$^*$}}
\newcommand{\sqdownend}{\multicolumn{1}{|c|}{$_+$}}
\newcommand{\squpdownend}{\multicolumn{1}{|c|}{$^*_+$}}

\newtheorem{theorem}{Theorem}
\newtheorem{lemma}{Lemma}[section]
\newtheorem{proposition}[lemma]{Proposition}
\newtheorem{corollary}[lemma]{Corollary}
\newtheorem{example}{Example}[section]
\newtheorem{definition}[example]{Definition}
\newtheorem{remark}[example]{Remark}
\newtheorem*{remark*}{Remark}
\newtheorem*{definition*}{Definition}

\begin{document}
\maketitle

\begin{abstract}
We give a simplification of Belkale's geometric proof of the Horn
conjecture.  Our approach uses the geometry of two-step
flag manifolds to explain the occurrence of the Horn inequalities
in a very straightforward way.  The arguments for both necessity and
sufficiency of the Horn inequalities are fairly 
conceptual when viewed in this framework.  We provide examples to
illustrate the method of proof.
\end{abstract}

\section{Introduction}
\subsection{General approach}

Horn's conjecture \cite{H} was originally formulated as a recursive 
method for solving a problem concerning the eigenvalues of
Hermitian matrices.  However, as a consequence of work of Klyachko
\cite{Kly}, Horn's conjecture can be 
reformulated as saying that the non-vanishing Schubert intersection 
numbers for Grassmannians satisfy a certain recursion.   This recursion
was first proved by Knutson and Tao \cite{KT}, using combinatorial 
methods.  Later, a geometric proof was given by 
Belkale \cite{B}, which was the inspiration for this paper.

For our purposes, a Schubert problem will refer to a collection
of Schubert varieties $\bar{\Omega}_{\sigma_1}, \ldots, 
\bar{\Omega}_{\sigma_s}$,
on some partial flag variety $G/P$.  A Schubert problem is non-vanishing
if the product of the corresponding cohomology classes is non-zero. 
Equivalently, a Schubert problem is non-vanishing if and only
if the general translates $\bar{\Omega}^{F_1}_{\sigma_1}, \ldots, 
\bar{\Omega}^{F_s}_{\sigma_s}$ of these Schubert varieties have
a non-empty generically transverse intersection.  This observation
allows one to study the vanishing question inside the tangent
space to the $G/P$.


In an effort to better understand the geometry behind Horn's conjecture,
we study the tangent spaces of Schubert 
varieties of two-step flag manifolds, and show how these are related
to the problem.  
The two-step flag manifold 
$$Fl(d,r,\CC^n) = \{ (S,V)\ |\ S \subset V \subset \CC^n,
\ \dim S = d,\ dim V = r\}$$
fibres over the Grassmannian $Gr(r,\CC^n)$, with fibre $Gr(d,V)$ at the
a point $V \in Gr(r,\CC^n)$.  A Schubert problem 
on $Gr(r,\CC^n)$ can be ``lifted'' in a number of 
different ways to a Schubert problem on $Fl(d,r,\CC^n)$, in such a way
that a non-empty transverse intersection on $Gr(r,\CC^n)$ lifts to 
a non-empty transverse
intersection on $Fl(d,r,\CC^n)$.  Each way of lifting corresponds to a 
non-vanishing Schubert problem inside the fibre $Gr(d,V)$.  However,
it is sometimes possible to see that the intersection on $Fl(d,r,\CC^n)$
is non-transverse for some trivial reasons.  These trivial conditions
are seen to be the Horn inequalities.

The difficult part of the Horn conjecture is to show the sufficiency
of the Horn inequalities.  In our approach this amounts to showing that 
a non-transverse intersection on
$Gr(r,\CC^n)$ lifts to something upstairs which is not only non-transverse,
but non-transverse for the aforementioned trivial reasons.  However,
once we have all the appropriate machinery in place, this turns 
out to be almost as straightforward as the ``easy'' direction of Horn's 
conjecture.

The two-step flag manifolds are not a necessary part
of the argument.  In principle the entire argument could be formulated
inside the tangent space of the Grassmannian.  This would probably
even lead to a shorter proof.  However, it is our opinion that the
geometry of the two-step flag manifold plays a fundamental role in
this picture, and to undermine its role would be remiss.

Although there are a number of new ideas in this paper, we do not
claim to be presenting an original or independent proof of the 
Horn conjecture.  Our objective in writing this was to better understand
the argument in \cite{B}, to simplify it in places, and to show
how it relates to some of our own previous work \cite{P1}.
A few of the results have been taken directly from \cite{B},
while some of the others merely contain old ideas which have been
dressed up in a new context.  We will try to indicate whenever possible
which ideas and results have been borrowed.

\subsection{Partial flag manifolds}

\subsubsection{Schubert varieties and Schubert positions}

Let $G = GL(n)$, and let $P \subset G$ be a parabolic subgroup, the
stabiliser subgroup of some $k$-step flag 
$$V^0 = \zeroset=V^0_0 \subsetneq V^0_1 \subsetneq \cdots 
\subsetneq V^0_k \subsetneq V^0_{k+1} = \CC^n.$$
We will assume for later convenience that the $V^0_i$ are coordinate
subspaces.
Let $d_j = \dim V_j$.
Then $G/P$ is a partial
flag manifold (of type $0<d_1<\cdots<d_k<n$).  

Let 
$$F = \zeroset=F_0 \subsetneq F_1 \subsetneq \cdots \subsetneq F_{n-1} 
\subsetneq F_n = \CC^n.$$
be a full flag on $C^n$. Let $B(F) \subset G$ denote its 
stabiliser.

Let $\sigma = \sigma_1 \ldots \sigma_n$ be a string of length $n$, where
$\sigma_l \in \{0, \ldots, k\}$, and the number of $j$'s in this string
is $d_{j+1}-d_j$.  It is a basic fact that $G/P$ has finitely many
$B(F)$-orbits, and that these orbits are indexed by the set of possible
$\sigma$ as follows.

Define 
$$\Omega^F_\sigma = \{V \in G/P\ \big|\ 
\sum_{j=1}^k \dim V_j \cap F_l - \dim V_j \cap F_{l-1} = \sigma_l\}.$$
$\Omega^F_\sigma$ is a $B(F)$-orbit, called the {\bf Schubert cell}
of $\sigma$ with respect to the flag $F$.  When the flag $F$ is
understood, or irrelevant, we shall also denote this $\Omega_\sigma$.
The {\bf Schubert variety} of $\sigma$, with respect to the flag $F$, is its
closure $\bar{\Omega}^F_\sigma$, which represents the Schubert class
$S^\sigma \in H^*(G/P)$.

\begin{example}
If the $\sigma_l$ are weakly decreasing, then $\Omega_\sigma$ consists of 
a single point which is a subflag of $F$.
If the $\sigma_i$ are weakly increasing, then $\Omega_\sigma$ is a dense
open subvariety of $G/P$, and $S^\sigma = 1 \in H^*(G/P)$.
\end{example}

In general we can easily calculate the dimension of this orbit:
$$\dim \Omega_\sigma = \#\{l<l'\ |\ \sigma_l < \sigma_{l'}\}.$$

If $V \in G/P$, then the {\bf Schubert position} of $V$ with respect
to the flag $F$ is the unique $\sigma$ such that 
$V \in \Omega^F_\sigma$.  When we have multiple flags $F^1, \ldots, F^s$
on $\CC^n$, the Schubert position of $V$ will be 
the $s$-tuple $(\sigma^1, \ldots \sigma^s)$, where $\sigma^i$ is the
Schubert position of $V$ with respect to $F^i$.

\subsubsection{$01$-strings from $\sigma$}
From $\sigma$ we construct {\em $01$-strings} in two different ways.
First, for each pair $(u,v)$, with $0 \leq u < v \leq k$, we define
a string $\sigma(uv)$.  This is constructed as follows: we
delete from $\sigma = \sigma_1 \ldots \sigma_n$ every
$\sigma_l \notin \{u,v\}$.  What remains will be a string consisting
only of `$u$'s and `$v$'s.  To convert this into a $01$-string,
we change every `$u$' to a \zero, and every `$v$' to a \one.

\begin{example}
If $\sigma=01312230132$, then $\sigma(13)$ is produced as follows:
$$\sigma = 0\mathbf{131}22\mathbf{3}0\mathbf{13}2
\mapsto \mathbf{131313} \mapsto 010101 = \sigma(13).$$
\end{example}

Each $\sigma(uv)$ corresponds to a Schubert cell on a different
Grassmannian.  In section \ref{sec:twostep} we investigate the
significance of these $\sigma(uv)$ on a two-step flag manifold.

In a completely different spirit, we define
strings $\sigma[j] = \sigma[j]_1 \ldots \sigma[j]_n$, 
for $1 \leq j \leq k$,
defined as follows.
$$ \sigma[j]_l = 
\begin{cases}
0,& \text{if $\sigma_l \leq k-j$} \\
1,& \text{if $\sigma_l > k-j$.}
\end{cases}
$$

These have a very natural geometric significance.  The
$k$-step flag variety has $k$ different projections onto 
Grassmannians $Gr(d_j, \CC^n)$.  The image of the Schubert
cell $\Omega_\sigma$ is projected onto $Gr(d_j, \CC^n)$
is the Schubert cell $\Omega_{\sigma[j]}$.

\begin{example}
If $\sigma=2103210$, then 
\begin{align*}
\sigma[1]&=0001000 \\
\sigma[2]&=1001100 \\
\sigma[3]&=1101110
\end{align*}
\end{example}

\subsubsection{Notation for Grassmannians}
When $P$ is a maximal parabolic, $G/P \cong Gr(r,n)$ 
is a Grassmannian,
and we can equivalently index the Schubert varieties and Schubert classes 
by partitions.  Let $\Lambda(r,n-r)$ denote the set of partitions
with $r$-parts, whose largest part is at most $n-r$, i.e.
$$\Lambda(r,n-r) = 
\{ \lambda = \mathbf{0} \leq \lambda_1 \leq \cdots \leq \lambda_r 
\leq \mathbf{n-r} \}.$$
There is a simple correspondence between these partitions and $01$-strings:
given a partition $\lambda \in \Lambda(r,n-r)$, one can construct a 
string as follows.

\begin{verse}
0. Begin with the empty string and $l=1$. \\
1. Append $\lambda_l - \lambda_{l-1}$ \zeros. \\
2. If $l=n-r$ stop. \\
3. Append a \one. \\
4. Increment $l$ and repeat steps 1-4.\\
\end{verse}

\begin{example}
$\lambda = \mathbf{0} \leq 0 \leq 1 \leq 3 \leq 3 \leq \mathbf{5}$, 
corresponds to the string $101001100$.
\end{example}

In the case of Grassmannians, we will sometimes find it more convenient to
use $\Lambda(r,n-r)$
to index our Schubert classes and Schubert varieties. We
therefore denote these $S^\lambda$ and $\bar{\Omega}_\lambda$ respectively.
It is worth noting at this time that 
$$\dim \Omega_\lambda = \sum_{l=1}^{n-r} \lambda_l =: |\lambda|.$$

\subsubsection{Induced flags}
Whenever we have a full flag $F$ on a vector space $W$, and $V \subset W$
is a subspace, we get induced flags $F_V$ and $F_{W/V}$ on $V$ and $W/V$
respectively.  These can be thought of as follows.

Consider the chain of subspaces
$$F_V = \bigg(
\zeroset=F_0 \cap V \subset F_1 \cap V \subset \cdots \subset F_{n-1} \cap V
\subset F_n \cap V = V \bigg).$$
Since $\dim V$ may be less than or equal to $\dim W$, this will no
longer be a flag, as some of the $F_{i-1} \cap V \subset F_i \cap V$
may become equalities.  However by eliminating all repeated elements,
one gets a full flag on $V$.  It is easy to check that the elements
we keep correspond precisely to the \ones in the Schubert position
of $V$ with respect to $F$.

Similarly we can construct a full flag on $W/V$.  Let $\Pi:W \to W/V$
denote the quotient map.  Then
$$F_{W/V}=
\bigg(
\zeroset=\Pi(F_0) \subset \Pi(F_1) \subset \cdots \subset \Pi(F_{n-1})
\subset \Pi(F_n) = W/V \bigg).$$
Again by eliminating repeated elements we obtain a full flag.
In this case, we keep the elements corresponding to the \zeros in the
Schubert position of $V$.

\subsection{The generalised Horn conjecture}
We can now state Horn's conjecture in the language of Schubert calculus.
Actually we will give the slightly more general statement which appears
in \cite{B}.

To simplify the notation in the statement a little, any time we
write $S^\lambda \in H^*(Gr(a,b))$ (or otherwise assume that $\lambda$ 
indexes such a class), then we implicitly assume 
$$\lambda = \mathbf{0} \leq \lambda_1  \leq \cdots \leq \lambda_a \leq
\mathbf{b-a} \in \Lambda(a,b-a).$$
 
\begin{definition}
\label{def:hornineq}
Let $\lambda^1, \ldots, \lambda^s \in \Lambda(r,n-r)$
Suppose that for every $d$, $1 \leq d \leq r$, 
and every non-zero product
$S^{\mu^1} \cdots S^{\mu^s} \in H^*(Gr(d,r))$, the inequality
\begin{equation}
\label{eqn:hornineq}
\sum_{i=1}^s \sum_{k=1}^d \lambda^i_{\mu^i_k +k} \geq (s-1)d(n-r)
\end{equation} 
holds.  In this case we say that the {\bf Horn inequalities} hold
for $\lambda^1, \ldots, \lambda^s$.
\end{definition}

\begin{theorem}[Generalised Horn conjecture]
The product 
$S^{\lambda^1} \cdots S^{\lambda^s} \in H^*(Gr(r,n))$ is non-zero
if and only if the Horn inequalities hold for
$\lambda^1, \ldots, \lambda^s$.
\end{theorem}

The standard form
of Horn's conjecture assumes that the product of 
$S^{\lambda^1} \cdots S^{\lambda^s}$ is of top degree, whereas
this formulation (which appears in \cite{B}) does not.  To accommodate this,
definition \ref{def:hornineq} allows an inequality for all
non-zero products in $H^*(Gr(d,r))$ rather than those just those
of top degree. This generalisation is easily shown to
imply the standard statement.

Definition \ref{def:hornineq}
differs slightly from the usual definition of the Horn inequalities
in another
small way.  Normally one does not include an inequality for the case 
$d=r$; one
simply uses $d$ for which $1 \leq d \leq r-1$.  The $d=r$ inequality
simply amounts to saying that 
$$\sum_{j=1}^s \codim \Omega_{\lambda^i} \leq \dim Gr(r,n),$$
i.e. this covers the case where the product vanishes for dimensional
reasons.  If we assume that $S^{\lambda^1} \cdots S^{\lambda^s}$ is of 
top degree then this inequality is always satisfied.

\section{Schubert calculus in the tangent space}
\subsection{General statements for $G/P$}
Let $F^1, \ldots, F^s$ be flags on $\CC^n$.
To determine whether
a product $S^{\sigma^1} \cdots S^{\sigma^s} \in H^*(Gr(r,n))$
vanishes, it is sufficient to consider whether the Schubert varieties
$\bar{\Omega}^{F^i}_{\sigma}$ have a point of intersection
when the flags $F^i$ are sufficiently generic.  This is due to 
the Kleiman-Bertini theorem \cite{Kl}, which tells us that 
if $F^i$ are generic, these Schubert varieties will intersect 
transversely (in positively oriented points).  Moreover, this point 
of intersection can be assumed
to be inside the open cell $\Omega^{F^i}_{\sigma}$

Thus we can take the following approach.  Choose the flags $F^i$
such that $V \in G/P$ is an intersection point of the Schubert
varieties $\Omega^{F^i}_{\sigma}$.  Subject to this restriction,
the flags $F^i$ should be as generic as possible.  We call generic
flags with this restriction {\bf almost generic} for $V$.  
The question then
becomes whether or not the Schubert varieties with respect to
almost generic flags intersect transversely.  
If they
intersect transversely, then 
Schubert varieties with respect to
fully generic flags must have a point of
intersection; however if they do not, then the point of intersection
is artificial, and fully generic flags will not give any
point of intersection.  

More importantly, this calculation can be
done on the level of tangent spaces.  We have the following basic
result which appears in $\cite{B, P1}$.

\begin{lemma}
\label{lem:tangent1}
$S^{\sigma^1} \cdots S^{\sigma^s} = 0 \in H^*(G/P)$ if and only
if the intersection of subspaces
$$\bigcap_{i=1}^s T_V \Omega^{F^i}_{\sigma^i} \subset T_V G/P$$
is non-transverse for $F^i$ almost generic.
\end{lemma}

Our $V$ can be any point of $G/P$, so let us take $V=V^0$.  Then
$P$ acts transitively on the flags $F^i$ such that 
$V^0 \in \Omega^{F^i}_{\sigma^i}$.  Thus to calculate 
$T_{V^0} \Omega^{F^i}_{\sigma^i}$, for generic $F^i$, it suffices
to compute it for a special $F^i$ and consider the action of a 
generic element of $P$ on the 
tangent space $T_{V^0} G/P = \frg/\frp$.  

There is always at least one
coordinate flag $\hat{F}(\sigma^i)$ such that 
$V^0 \in \Omega^{\hat{F}(\sigma^i)}_{\sigma^i}$; we will use one of these.
The standard flag on $\CC^n$ is
$$\Fstand = \zeroset \subsetneq \langle x_1 \rangle
\subsetneq \langle x_1, x_2 \rangle \subsetneq \cdots 
\subsetneq \langle x_1, x_2, \ldots, x_{n-1} \rangle $$
The Weyl group $S_n$ acts transitively on the set of coordinate flags 
of $\CC^n$, thus $\hat{F}^i = w_i \cdot \Fstand$ for some $w_i \in S_n$.
Of all possible choices for the $\hat{F}(\sigma^i)$, we will always
choose $\hat{F}(\sigma^i) = w_i \cdot \Fstand$ with $w_i$ minimal
in the Bruhat order.

Denote the tangent space to the Schubert variety 
$\Omega^{\hat{F}(\sigma)}_{\sigma}$ at $V_0$ by
$$\hat{Z}_\sigma = T_{V^0} \Omega^{\hat{F}(\sigma)}_{\sigma} 
\subset \frg/\frp.$$
Let $Z_\sigma$ denote a generic $P$-translate of $\hat{Z}_\sigma$.  This
notation will be convenient, as we will often need to consider 
intesections
$\bigcap_{i=1}^s p_i \hat{Z}_{\sigma_i}$ where $p_i \in P$ are 
generic; we can now write this simply as $\bigcap_{i=1}^s Z_{\sigma_i}$.
For the special cases of Grassmannians and two-step flag manifolds, we will
use the letters $X$ or $Y$ respectively instead of $Z$.

Lemma \ref{lem:tangent1} can be stated equivalently as follows:

\begin{lemma}
\label{lem:tangent2}
$S^{\sigma^1} \cdots S^{\sigma^s} = 0 \in H^*(G/P)$ if and only
if the intersection of subspaces
$$\bigcap_{i=1}^s Z_{\sigma^i} \subset \frg/\frp$$
is non-transverse.
\end{lemma}

\begin{remark}
Although we have eliminated the flags $F^i$ from this statement,
we will sometimes wish to think of $Z_{\sigma^i}$ as being determined
by generic flags, rather than as generic translates of 
$\hat{Z}_{\sigma^i}$.
\end{remark}

\begin{remark}
Whenever we use the notation $Z_{\sigma^i}$, we tacitly assume (unless
otherwise specified) that underlying flags $F^i$ are almost generic.
If the $F^i$ are almost generic, we say that
$Z_{\sigma^1}, \ldots, Z_{\sigma^s}$ are in {\bf general position}.
This is, of course, equivalent to the statement that
$Z_{\sigma^i}$ is a \emph{generic} $P$-translate of 
$\hat{Z}_{\sigma^i}$.
\end{remark}

To calculate a particular $\hat{Z}_{\sigma}$, we use the following fact.

\begin{proposition}
\label{prop:calcZ}
$\hat{Z}_{\sigma}$ is the image of $\frb(\hat{F}(\sigma))$ under
$\frg \to \frg/\frp$.
\end{proposition}

\begin{proof}
As $\Omega^{\hat{F}(\sigma)}_{\sigma}$ is a $B(\hat{F}(\sigma))$ orbit,
the tangent space is generated by $\frb(\hat{F}(\sigma))$.
\end{proof}

\subsection{Grassmannians}
\label{sec:grassmannian}

\subsubsection{The tangent space to a Grassmannian Schubert variety}

To distinguish the
Grassmannian as a special case, we use the notation
$\hat{X}_{\lambda}$ for $\hat{Z}_\sigma$, 
and $X_\lambda$ for the generic translate $X_\sigma$,
where $\lambda$ is the partition corresponding to the $01$-string $\sigma$.
We will take $V^0$ to be the coordinate subspace 
$\langle x_{n-r+1}, \ldots, x_n \rangle \in Gr(r,n)$. 

We will now identify the subspace $\hat{X}_\lambda$.
Put $V = V^0$, and $Q = \CC^n/V$.
Now $\frg/\frp$ can be naturally identified with $\Hom(V, Q)$,
so we can view $\hat{X}_\lambda$ as a set of homomorphisms $\phi:V \to Q$.
Both $V$ and $Q$ inherit full flags from the standard flag.  $V$ inherits the
flag 
$$\Fstand_V = \zeroset \subsetneq \langle x_{n-r+1} \rangle
\subsetneq \langle x_{n-r+1}, x_{n-r+2} \rangle \subsetneq \cdots
\subsetneq \langle x_{n-r+1}, \ldots, x_{n-1} \rangle \subsetneq V$$
and $Q$ inherits the image of
$$\Fstand_Q = V \subsetneq V + \langle x_{1} \rangle
\subsetneq V + \langle x_{1}, x_{2} \rangle \subsetneq \cdots
\subsetneq V + \langle x_{1}, \ldots, x_{n-r-1} \rangle \subsetneq \CC^n$$
under the quotient map $\CC^n \to Q$.

\begin{proposition}
Under these identifications,
$$\hat{X}_\lambda = 
\{\phi \in \Hom(V,Q)\ |\ \phi((\Fstand_V)_l) \subset
(\Fstand_Q)_{\lambda_l}\ l=1, \ldots, r\}.
$$
\end{proposition}

\begin{proof}
Let $\sigma$ denote the $01$-string corresponding to $\lambda$, and put
$$
\alpha_k = 
\begin{cases}
 \sum_{l=1}^k 1-\sigma_l, &\text{if $\sigma_k = 0$} \\
n-r+\sum_{l=1}^k \sigma_l, &\text{ if $\sigma_k = 1$.} 
\end{cases}
$$
So $\alpha_k \leq n-r$ if $\sigma_k=0$, and $\alpha_k >n-r$ if $\sigma_k=1$.
Then the flag $\hat{F}(\sigma)$ is 
$$\hat{F}(\sigma)
= \zeroset \subsetneq \langle x_{\alpha_1} \rangle
\subsetneq \langle x_{\alpha_1}, x_{\alpha_2}\rangle \subsetneq \cdots
\subsetneq \langle x_{\alpha_1}, \ldots, x_{\alpha_{r-1}} \rangle
\subsetneq \CC^n.$$
To see this, note that $V \cap \hat{F}(\sigma)_l$ jumps in dimension
exactly when $\sigma_l=1$, and $\alpha$ is the smallest permutation that
accomplishes this.

By proposition \ref{prop:calcZ}, $\hat{X}_\sigma$ is is identified with
$\frb(\hat{F}(\sigma))/\frp$.
Let $\frn$ denote the orthogonal complement to $\frp$.
Since $V$ and $\hat{F}(\sigma)$ are coordinate flags, we have
$\frb(\hat{F}(\sigma))/\frp = \frb(\hat{F}(\sigma)) \cap \frn/\frp$.
Thus it suffices to determine $\frb(\hat{F}(\sigma)) \cap \frn$.
We see that an element $\phi \in \frn$ preserves $\hat{F}(\sigma)$
if and only if 
for $j \in \{n-r+1, \ldots, n\}$, 
$$\phi(x_j) \subset \langle x_{\alpha_1}, x_{\alpha_2}, \ldots, 
x_{\alpha_j} \rangle$$
But since 
$\Image(\phi) \subset \langle x_1, \ldots, x_{n-r} \rangle$, this
is equivalent to 

$$\phi(x_j) \subset V + \langle x_1, \ldots, x_l \rangle
= (\Fstand_Q)_l$$
where $l$ is the number of \zeros before the $j^{\rm th}$ \one.  Thus
$l = \lambda_j$.
\end{proof}

\begin{example}
In $Gr(4, \CC^9)$,
if $\lambda = \mathbf{0} \leq 0 \leq 1 \leq 3 \leq 3 \leq \mathbf{5}$,
then with respect to the standard coordinate bases on $V$ and $Q$,
$$\hat{X}_\lambda = \
\left (
\begin{matrix}
0 & * & * & * \\
0 & 0 & * & * \\
0 & 0 & * & * \\
0 & 0 & 0 & 0 \\
0 & 0 & 0 & 0 
\end{matrix}
\right)
$$
In general, if we write $\hat{X}_\lambda$ in this form, the number of
$*$'s in column $j$ will be $\lambda_j$.
\end{example}

The action of $P$ on $\Hom(V,Q)$ is given by
$$
\left(
\begin{matrix}
A & 0 \\
C & B 
\end{matrix}
\right) 
\cdot \phi
= A \phi B^{-1}.
$$
To compute the generic intersection of two subspaces 
$X_\lambda \cap X_\mu$, it suffices to take $X_\lambda = \hat{X}_\lambda$, 
and $X_\mu = p \cdot \hat{X}_\mu$, where
$$ p =
\left(
\begin{matrix}
 &       &1& & &   \\
 & \revddots & & &\mbox{\Huge $0$}&  \\
1&       & & & &   \\
 & &       & &       &1   \\
 &\mbox{\Huge $0$}&       & &\revddots&   \\
 & &       &1&       &    \\
\end{matrix}
\right)
$$
This is the tangent space analogue of intersecting a Schubert variety
with an opposite Schubert variety.

\begin{example}
\label{ex:oppositeintersection}
If $\lambda = \mathbf{0} \leq 0 \leq 1 \leq 3 \leq 3 \leq \mathbf{5}$,
$\mu = \mathbf{0} \leq 3 \leq 3 \leq 3 \leq 5 \leq \mathbf{5}$, then
\begin{align*}
X_\lambda \cap X_\mu &=
\left (
\begin{matrix}
0 & * & * & * \\
0 & 0 & * & * \\
0 & 0 & * & * \\
0 & 0 & 0 & 0 \\
0 & 0 & 0 & 0 
\end{matrix}
\right)
\cap
\left (
\begin{matrix}
* & 0 & 0 & 0 \\
* & 0 & 0 & 0 \\
* & * & * & * \\
* & * & * & * \\
* & * & * & * 
\end{matrix}
\right) \\
&= 
\left (
\begin{matrix}
0 & 0 & 0 & 0 \\
0 & 0 & 0 & 0 \\
0 & 0 & * & * \\
0 & 0 & 0 & 0 \\
0 & 0 & 0 & 0 
\end{matrix}
\right)
\end{align*}
The intersection is codimension $18$; however the expected codimension
is $\codim X_\lambda + \codim X_\mu = 13+6 =19$.  Thus 
$S^\lambda S^\mu = 0 \in H^*(Gr(r,n))$.
\end{example}

\subsubsection{$X_\lambda$ as a set of homomorphisms}
\label{sec:grassmannianhoms}
Let $\phi \in X_\lambda \subset \Hom(V,Q)$.  
This represents a direction in which we can
perturb $V$ and remain in the tangent space to the Schubert
variety.  The kernel, $\ker \phi$ represents the maximal subspace of $V$ which
is preserved by this perturbation.  One of the key ideas which we take
from Belkale's proof is to examine $\ker \phi$ for a generic 
$\phi \in \bigcap_{i=1}^s X_{\lambda^i}$.  For now, however, let us 
restrict our attention to a single $X_\lambda$.

Let us recall that the construction of $X_\lambda$ was actually a
tangent space of a Schubert variety relative
to some generic flag $F$ on $\CC^n$.  Thus $V$ carries an induced
flag $F_V$, and $Q$ carries an induced flag $F_Q$.  The action
of $GL(V) \times GL(Q)$ on $\Hom(V,Q)$ corresponds to changing
these induced flags.  $F_V$ and $F_Q$ carry all the
relevant information (for our purposes) about the original flag
$F$, and moreover, the action of $GL(V) \times GL(Q)$ on 
$Flags(V) \times Flags(Q)$ is transitive.  Thus we shall stop
thinking about $F_V$ and $F_Q$ as the induced flags of $F$, and
instead think of them as an independent pair of generic flags
on $V$ and $Q$ respectively.

\begin{remark}
\label{rmk:genericinducedflags}
We have in fact just observed that the induced flags on $V$
and $\CC^n/V$ are generic, if $V \in Gr(r,n)$ is a generic
point of intersection of Schubert varieties in general position.
\end{remark}

Let $S$ be a subspace of $V$ whose Schubert position relative to the
flag $F_V$ is $\rho$.  Put $d = \dim S$. 
If $\phi \in \Hom(V,Q)$ and $\ker \phi
\supset S$, then $\phi$ descends to a map $\tilde{\phi} \in \Hom(V/S,Q)$.
We now consider the space 
$$X_\lambda\modmod S = \{\tilde\phi \in \Hom(V/S,Q)\ |\ \phi \in X_\lambda,
\ \ker\phi \supset S\} = X_\lambda \cap \Hom(V/S,Q).$$

\begin{lemma}
$X_\lambda\modmod S = X_{\lambda'}$, where $\lambda' \in \Lambda(r-d, n-r)$
is the subpartition
$$\lambda' = \mathbf{0} \leq \lambda_{k_1} \leq \cdots \leq 
\lambda_{k_{r-d}} \leq \mathbf{n-r}$$
where $k_j$ are the positions of the \zeros of $\rho$.
The flag on $V/S$ associated to $X_{\lambda'}$ is the induced
flag $(F_V)_{V/S}=:F_{V/S}$.
\end{lemma}

\begin{proof}
The action of $P$ simply gives a change of basis on $V$
and $Q$.  Thus it suffices to prove this in the case where 
$F = \hat{F}(\lambda)$, in which case $X_\lambda = \hat{X_\lambda}$ and
$F_V = \Fstand_V$.

We must verify that for every $\phi \in X_\lambda\modmod S$, that 
$\phi ((F_{V/S})_j) \subset (F_Q)_{\lambda'_j}$.  However, this is
straightforward, as $(F_{V/S})_j$ is the image under the quotient map
of $(F_V)_{k_j}$ which maps to 
$(F_Q)_{\lambda_{k_j}} =(F_Q)_{\lambda'_j}$.  Thus $X_\lambda\modmod S \subset
X_{\lambda'}$.  

Moreover, every element
of $X_{\lambda'}$ can be seen to have a unique lifting to 
$\{\phi \in X_\lambda\ |\ S \supset \ker \phi\}$.  To see this, note
that by a change of coordinates which preserves $F_V$ we can make
$S$ is a coordinate subspace, in which case, this is 
obvious.  Thus $X_\lambda\modmod S = X_{\lambda'}$.
\end{proof}

\subsubsection{Genericity of $S$}
The one remaining fact we will need is that when $S$ is the special
subspace $S = \ker \phi$, $\phi \in \bigcap_{i=1}^s X_{\lambda^i}$, 
then
the flags $F_S$ and  $F_{V/S}$ are actually {\em generic} flags.  
It is a priori
conceivable that by choosing this particular $S$ (its definition 
involves the flags $F^i$), we could have undone all the genericity
we had before, thereby
ending up in a position where all the flags $F^i_{V/S}$ are not
generic with respect to each other (for example, they could all be equal).
This would be most unfortunate; however luckily it does not happen.

Let $\rho^i$ denote the Schubert position of $S \subset V$ with
respect to the flag $F^i_V$.  Thus 
$S \in \bigcap_{i=1}^s \Omega^{F^i_V}_{\rho^i}$.  We must show 
that $S$ is in fact a generic point of this intersection of Schubert
varieties.  This is sufficient, since the induced flags at
a generic point of an intersection of Schubert varieties are generic
(c.f. remark \ref{rmk:genericinducedflags}).

The idea is to fix generic flags $F^i_V$, 
while the flags $F^i_Q$ vary.
We show that there cannot exist a subvariety 
$T \subset \bigcap_{i=1}^s \Omega^{F^i_V}_{\rho^i}$ such that
$\ker \phi \in T$, for a generic choice of 
$\phi \in \bigcap_{i=1}^s X_{\lambda^i}$, and generic induced flags 
$F^i_Q$.  This is a consequence of a generalisation of the Kleiman
moving lemma, due to Belkale.

\begin{lemma}[{\rm Belkale \cite{B}}]
\label{lem:moving}
Let $H$ be an algebraic group acting on a variety $X$.
Suppose $\pi: X \to Y$ is an $H$-invariant fibration, such that $H$
acts transitively on the fibres.  Let $Z_i \subset X$, and 
$Y_i=\pi(Z_i) \subset Y$ $i \in \{1, \ldots, s\}$, be subvarieties,
such that $\pi|_{Z_i}: Z_i \to Y_i$ is a fibration.  Put 
$Y_0 = \bigcap_{i=1}^s Y_i$, and $X_0 = \pi^{-1}(Y_0)$.  Let 
$T \subset Y_0$, be a subvariety.  Then for generic
$h_i \in H$, the intersection
$$\pi^{-1}(T) \cap \left ( \bigcap_{i=1}^s h_i \cdot Z_i \right )$$
has the expected dimension as an intersection inside $X_0$.
That is,
\begin{equation}
\label{eqn:moving}
\codim \pi^{-1}(T) \cap \left ( \bigcap_{i=1}^s h_i \cdot Z_i \right )
= \codim \pi^{-1}(T) + 
\sum_{i=1}^s \codim (Z_i \cap X_0)
\end{equation}
(here $\codim$ means ``codimension inside $X_0$'').
\end{lemma}

In this example, 
$$X = \Hom_d(V,Q) 
= \{\psi \in \Hom(V,Q)\ |\ \dim \ker \psi = d\},$$
$$Y = Gr(n-d,V), \qquad \pi(\psi) = \ker(\psi)$$
The group $H$ is $GL(Q)$, which we note preserves the kernel of
$\psi \in Z_i$.  The action of $GL(Q)$ acts transitively on the
induced flags $F^i(Q)$.

The subvariety $Z_i \subset X_{\lambda^i} \subset X$ consists of 
those elements of 
$X_{\lambda^i}$ whose kernel is in Schubert position $\rho^i$.
Thus a generic translate $h_i \cdot Z_i$ is simply $X_{\lambda^i}$
for a different generic choice of flags $F^i_Q$.
Thus the image $\pi(\psi)$ of a generic point in the intersection 
$\psi \in \bigcap_{i=1}^s h_i \cdot Z_i$, is just the kernel of a
generic element of $\bigcap_{i=1}^s X_{\lambda^i}$.

The image of $Z_i$ under $\pi$ is the open Schubert cell
$Y_i = \Omega^{F^i_V}_{\rho^i}$.  We wish to show that there is 
no subvariety $T \subset \bigcap Y_i$ such that 
$\ker \psi \in T$ for all choices above.  In other words, for
any subspace $T$ we can choose $\psi$ and generic flags $F^i_Q$ 
such that $\psi \in \bigcap_{i=1}^s X_{\lambda^i}$ and
$\ker \psi \notin T$.  Equivalently, we must show there exists 
$\psi \in \bigcap h_i \cdot Z_i$, with $\psi \notin \pi^{-1}(T)$.
But this is clear from lemma \ref{lem:moving}, since otherwise 
the codimensions in equation (\ref{eqn:moving}) would not add.

\begin{corollary}
\label{cor:genericdescent}
Let $S = \ker \phi$ for a generic $\phi \in \bigcap_{i=1}^s 
X_{\lambda^i}$.  Then the $X_{\lambda^i}\modmod S = X_{{\lambda^i}'}$ 
are in general position.  Moreover
there is an element $\tilde \phi \in \bigcap_{i=1}^s X_{{\lambda^i}'}
\subset \Hom(V/S,Q)$ such that $\ker \tilde \phi = \zeroset$.
\end{corollary}

\begin{proof}
The fact that the $X_{\lambda^i}\modmod S$ are generic simply means that
the induced flags $F^i_{V/S}$ and $F^i_{Q}$ are generic, which is
what we have just shown.  For the second statement, since 
$\ker \phi = S$, $\phi$ descends to well defined map
$\tilde \phi:V/S \to Q$ with 
$\tilde \phi \in \bigcap_{i=1}^s X_{\lambda^i}\modmod S$ and
$\ker \tilde \phi = \zeroset$.
\end{proof}

\begin{corollary}
\label{cor:transversereduction}
The intersection $\bigcap_{i=1}^s X_{\lambda^i}\modmod S$ is transverse.
\end{corollary}

The argument can be seen as special case of \cite{B}[lemma 2.18].

\begin{proof}
Since the $X_{\lambda^i}\modmod S = X_{{\lambda^i}'}$ 
are linear subspaces of $\Hom(V/S,Q)$ their
intersection is necessarily equidimensional.  Thus it suffices to
show that the intersection is transverse on an open subset
which contains a point of intersection.

Consider the space 
$\Hom_0(V/S,Q) = \{\psi \in \Hom(V/S,Q)\ |\ \ker \psi = \zeroset\}.$
This is a homogeneous space under the action of $GL(V/S) \times GL(Q)$,
and it is a Zariski open subset of $\Hom(V/S,Q)$.  
By corollary \ref{cor:genericdescent}
the $X_{{\lambda^i}'}$ are generic translates of $\hat X_{{\lambda^i}'}$
by elements of $GL(V/S) \times GL(Q)$.  So by Kleiman-Bertini, the
intersection 
$$\bigcap_{i=1}^s \left(\Hom_0(V/S,Q) \cap X_{\lambda^i}\right)$$
is a transverse intersection. Moreover, it contains the point 
$\tilde \phi$
which shows that the intersection $\bigcap_{i=1}^s X_{{\lambda^i}'}$
is transverse.
\end{proof}

\subsection{Two-step flag manifolds}
\label{sec:twostep}

\subsubsection{The tangent space to a Schubert variety in a two-step
flag manifold}

Let $0<d<r<n$.  We consider the two-step flag manifold $Fl(d,r,\CC^n)$.
For two-step flag manifolds, we use the notation
$\hat{Y}_{\sigma}$ for $\hat{Z}_\sigma$, 
and $Y_\sigma$ for the generic translate $Z_\sigma$.
Here $\sigma$ is a $012$-string, with $d$ \twos, $r-d$ \ones, and
$n-r$ \zeros.  Let $\eta_1 < \eta_2 < \cdots <\eta_{n-r}$
denote the positions of the \zeros (i.e $\sigma_{\eta_m}=0$, for
$m \leq n-r$).
Let $\eta_{n-r+1} < \cdots < \eta_{n-d}$ denote the positions
of the \ones, and $\eta_{n-d+1} < \cdots < \eta_{n}$ denote the 
positions of the \twos.

Our first objective is to describe the space $\hat{Y}_{\sigma}$.
Let $V$ denote the coordinate subspace
$V = \langle x_{n-r+1}, \ldots, x_n \rangle$, and $S$ denote the
coordinate subspace 
$S = \langle x_{n-d+1}, \ldots, x_n \rangle$.  (Eventually $d$ and
$S$ will play the same role as they did in section 
\ref{sec:grassmannian}.)  We take our base flag $V^0$ to be
the coordinate two-step flag
$$V^0 = \zeroset \subset S \subset V \subset \CC^n.$$

Now $\frg/\frp$ has a basis descending from the standard basis
$$\big\{ E_{jk}
\ \big|\ \mbox{$j \leq n-r$, $k>n-d$, and ($k>n-d$ or $j \leq n-r$)}
\big\}$$
where $E_{jk}$ is the image under $\frg \to \frg/\frp$ of the $n \times n$
matrix whose only non-zero entry is a \one in the $(j,k)$-position.
These basis vectors naturally partition $\frg/\frp$ into three blocks:
the upper left block is spanned by those $E_{jk}$ such that
$1 \leq j \leq n-r$, $n-r < k \leq n-d$; the lower right block is
spanned by those $E_{jk}$ such that 
$n-r < j \leq n-d$, $n-d < k \leq n$, and the upper right block is
spanned by those $E_{jk}$ such that 
$1 \leq j \leq n-r$, $n-d < k \leq n$.  The first two are 
naturally viewed as subspaces, while the last is more naturally
viewed as a quotient space.

Since $S$ and $V$ are coordinate subspaces, $\hat{Y}_{\sigma}$ is
spanned by some subset of the $E_{jk}$.

\begin{proposition}  
$\hat{Y}_{\sigma} = \{E_{jk}\ |\ \eta_j < \eta_k\}$.
\end{proposition}

\begin{proof}
$\eta = \eta_1 \ldots \eta_n$ is the element of the $S_n$ 
such that $\eta^{-1} \cdot \Fstand = \hat{F}(\sigma)$.  Thus
$\frb(\hat{F}(\sigma)) = \eta^{-1} \cdot \frb(\Fstand)$.  One 
can easily check that the image in $\frg/\frp$ is
$\{E_{jk}\ |\ \eta_j < \eta_k\}$.
\end{proof}

\begin{example}
\label{ex:Y}
Let $d=2$, $r=5$, $n=9$.  Then $\frg/\frp$ looks like
$$
\left(
\begin{matrix}
\cdot & \cdot & \cdot & \cdot & * & * & * & * & * \\
\cdot & \cdot & \cdot & \cdot & * & * & * & * & * \\
\cdot & \cdot & \cdot & \cdot & * & * & * & * & * \\
\cdot & \cdot & \cdot & \cdot & * & * & * & * & * \\
\cdot & \cdot & \cdot & \cdot & \cdot & \cdot & \cdot & * & * \\
\cdot & \cdot & \cdot & \cdot & \cdot & \cdot & \cdot & * & * \\
\cdot & \cdot & \cdot & \cdot & \cdot & \cdot & \cdot & * & * \\
\cdot & \cdot & \cdot & \cdot & \cdot & \cdot & \cdot & \cdot & \cdot \\
\cdot & \cdot & \cdot & \cdot & \cdot & \cdot & \cdot & \cdot & \cdot 
\end{matrix}
\right)
$$
Let $\sigma = 021010201$.  Then $\eta = 146835927$, and 
$$\hat{Y}_{\sigma} =
\left(
\begin{matrix}
\cdot & \cdot & \cdot & \cdot & * & * & * & * & * \\
\cdot & \cdot & \cdot & \cdot & 0 & * & * & 0 & * \\
\cdot & \cdot & \cdot & \cdot & 0 & 0 & * & 0 & * \\
\cdot & \cdot & \cdot & \cdot & 0 & 0 & * & 0 & 0 \\
\cdot & \cdot & \cdot & \cdot & \cdot & \cdot & \cdot & 0 & * \\
\cdot & \cdot & \cdot & \cdot & \cdot & \cdot & \cdot & 0 & * \\
\cdot & \cdot & \cdot & \cdot & \cdot & \cdot & \cdot & 0 & 0 \\
\cdot & \cdot & \cdot & \cdot & \cdot & \cdot & \cdot & \cdot & \cdot \\
\cdot & \cdot & \cdot & \cdot & \cdot & \cdot & \cdot & \cdot & \cdot 
\end{matrix}
\right )$$
Note that the above pictures of  $\frg/\frp$ and $\hat{Y}_\sigma$
break up into three blocks,
$$ \left(
\begin{matrix}
* & * & * \\
0 & * & * \\
0 & 0 & * \\
0 & 0 & * 
\end{matrix}
\right)
\qquad
\left(
\begin{matrix}
* & * \\
0 & * \\
0 & * \\
0 & 0
\end{matrix}
\right)
\qquad
\left(
\begin{matrix}
0 & * \\
0 & * \\
0 & 0
\end{matrix}
\right )
$$
and each block contains a Young-diagram shaped picture.  These
pictures are $\hat{X}_{\sigma(01)}$, $\hat{X}_{\sigma(02)}$ and
$\hat{X}_{\sigma(12)}$.
\end{example}

\subsubsection{Grassmannian problems in the two-step flag manifold}
In example \ref{ex:Y} one can see that $\frg/\frp$ consists of
three rectangular blocks, and the diagram representing 
$\hat{Y}_\sigma$ is shaped like a Young diagram when restricted to
each of these blocks.  Thus these are actually diagrams for
some $\hat{X}_\lambda$.

Inside $\frg/\frp$ there are two natural subspaces, corresponding
to the fibres of the forgetful maps $G/P \to Gr(r,\CC^n)$ 
and $G/P \to Gr(d,\CC^n)$.  The subspaces are $\Hom(S,V)$,
and $\Hom(V/S, Q)$ (where $Q = \CC^n/V$).  In terms of our basis 
for $\frg/\frp$, 
these correspond to the lower right
and upper left blocks respectively.  Note that these subspaces
are in fact $P$-invariant.  Thus the quotient map
$$\Pi_{\Hom(S,Q)}  : \frg/\frp \to \Hom(S, Q) =
(\frg/\frp)\big/ \big(\Hom(S,V) \oplus \Hom(V/S,Q) \big )$$
is $P$-equivariant.

By simply noting that the dimensions are correct we see that is
it possible to view $X_{\sigma(01)}$ as a subspace of $\Hom(V/S,Q)$,
$X_{\sigma(12)} \subset \Hom(S,V)$, and 
$X_{\sigma(02)} \subset \Hom(S,Q)$.  However, better than this,
we have the following result.

\begin{proposition} We have the following identifications.
\begin{enumerate}
\item $Y_\sigma \cap \Hom(V/S,Q) = X_{\sigma(01)}$.
\item $Y_\sigma \cap \Hom(S,V) = X_{\sigma(12)}$.
\item $\Pi_{\Hom(S,Q)}(Y_\sigma) = X_{\sigma(02)}$.
\end{enumerate}
Moreover, if $Y_{\sigma^i}$, $i \in \{1, \ldots, s\}$
are in general position, 
then so are $X_{\sigma^i(01)}$, $X_{\sigma^i(12)}$,
and $X_{\sigma^i(02)}$.
\end{proposition}
\begin{proof}
Since the spaces $\Hom(V/S,Q), \Hom(S,V) \subset \frg/\frp$ are 
$P$-invariant, it suffices to check statements 1--3 for $\hat{Y}_\sigma$.

For statement 1, let $\zeta_1< \cdots < \zeta_{n-d}$, and 
$\zeta_{n-d+1} < \cdots <\zeta_{n-r}$ denote the
positions of the \zeros and \ones respectively for $\sigma(01)$, and
let $\lambda$ denote the partition corresponding to $\sigma(01)$.
Note that $j < \lambda_k$ if and only if
$\zeta_j < \zeta_k$ if and only if $\eta_j < \eta_k$.
Thus we have 
$$\hat{Y}_\sigma \cap \Hom(V/S,Q) =
\langle E_{jk}\ |\ 
\mbox{$1 \leq j \leq n-r$, $n-r < k< n-d$, and $\eta_j < \eta_k$} 
\rangle$$
On the other hand,
\begin{align*}
\hat{X}_{\sigma(01)} &=
\langle E_{jk}\ |\ E_{jk}((\Fstand_{V/S})_l) \subset 
(\Fstand_Q)_{\lambda_l},\ \forall l \rangle \\
&=
\langle E_{jk}\ |\ E_{jk}((\Fstand_{V/S})_k) \subset 
(\Fstand_Q)_{\lambda_k} \rangle \\
&=
\langle E_{jk}\ |\ (\Fstand_Q)_j \subset 
(\Fstand_Q)_{\lambda_k} \rangle \\
&=
\langle E_{jk}\ |\ j < \lambda_k \rangle  \\
&=
\langle E_{jk}\ |\ \zeta_j < \zeta_k \rangle  
\end{align*}  
To see the genericity, note that $P$ acts on $\Hom(V/S,Q)$,
and contains the subgroup $GL(V/S) \times GL(Q)$ (acting in the
usual way).  Thus if $Y_{\sigma^i}$ are generic, so are 
$X_{\sigma^i(01)}$.

A similar argument holds for $\sigma(12)$ and $\sigma(02)$. 
\end{proof}

Even more importantly, we can sometimes use the transversality
(or lack thereof) for these Grassmannian Schubert problems to deduce
transversality in the two-step flag manifold.

\begin{lemma}
\label{lem:splitting}
Let $\sigma^1, \ldots, \sigma^s$ be $012$-strings, giving rise
to Schubert classes on the two-step flag manifold $Fl(d,r,\CC^n)$.
\begin{enumerate}
\item
If $\bigcap_{i=1}^s X_{\sigma^i(02)}$ is non-transverse, then
$\bigcap_{i=1}^s Y_{\sigma^i}$ is non-transverse.
\item
If $\bigcap_{i=1}^s X_{\sigma^i(01)}$,
$\bigcap_{i=1}^s X_{\sigma^i(12)}$, and
 $\bigcap_{i=1}^s X_{\sigma^i(02)}$ are all transverse, then
$\bigcap_{i=1}^s Y_{\sigma^i}$ is transverse.
\end{enumerate}
\end{lemma}

\begin{proof}
These follow from elementary facts about intersections in quotient
spaces and subspaces.
\end{proof}

One important special case of this theorem is the following, which
is a variant on the vanishing criterion in \cite{P1}.
\begin{corollary}
\label{cor:almosthorn}
If $\sum_{i=1}^s \dim X_{\sigma^i(02)} < (s-1)d(n-r)$ then
$\bigcap_{i=1}^s Y_{\sigma^i}$ non-transverse.
\end{corollary}

\begin{proof}
If $\sum_{i=1}^s \dim X_{\sigma^i(02)} < (s-1)d(n-r)$ then 
$\bigcap_{i=1}^s X_{\sigma^i(02)}$ is non-transverse for dimensional
reasons.
\end{proof}

\subsubsection{The strings $\sigma[1]$ and $\sigma[2]$}
We can view $\Hom(V,Q)$ and $\Hom(S, \CC^n/S)$ as quotient
spaces of $\frg/\frp$ as well:  $\Hom(V,Q) = (\frg/\frp)\big/\Hom(S,V)$
and $\Hom(S, \CC^n/S) = (\frg/\frp) \big/\Hom(V/S,Q)$.
Let $\Pi_{\Hom(V,Q)}$ and $\Pi_{\Hom(S, \CC^n/S)}$ denote the 
projection maps.
Counting \ones and \zeros, we see $X_{\sigma[2]}$ and $X_{\sigma[1]}$
are of the correct dimensions to view
$X_{\sigma[2]} \subset \Hom(V,Q)$ and 
$X_{\sigma[1]} \subset \Hom(S, \CC^n/S)$.

\begin{proposition}
$\Pi_{\Hom(V,Q)}(Y_\sigma)$ is an $X_{\sigma[2]}$, and
$\Pi_{\Hom(S,\CC^n)}(Y_\sigma)$ is an $X_{\sigma[1]}$; however,
they are not in general position.
\end{proposition}

\begin{proof}
Let $F$ be a flag such that $(S,V) \in \Omega^F_{\sigma}$.  Then
$V \in \Omega^F_{\sigma[2]}$, and $S \in \Omega^F_{\sigma[1]}$.
It follows that $Y_\sigma = \frb(F)/\frp$ maps surjectively
onto $X_{\sigma[2]} = \frb(F)/\Stab(V)$ and 
$X_{\sigma[1]} = \frb(F)/\Stab(S)$.  The quotient maps are
precisely $\Pi_{\Hom(V,Q)}$ and $\Pi_{\Hom(S, \CC^n/S)}$.

They are not in general position, since the flag $F$ is
almost generic for the two step flag 
$\zeroset \subset S \subset V \subset \CC^n$, but not
necessarily almost generic for $V$ or $S$.
\end{proof}

We can now tie this to the material in section 
\ref{sec:grassmannianhoms}.
\begin{proposition}  
Fix a flag $F$ such that $(S,V) \in \Omega^{F}_{\sigma}$.
View $X_{\sigma[2]} \subset \Hom(V,Q)$ and $X_{\sigma(01)} \subset \Hom(V/S,Q)$,
as being tangent spaces to Schubert varieties relative the flag $F$
(or its induced flags).  Then
$X_{\sigma[2]}\modmod S = X_{\sigma(01)}$.
\end{proposition}

\begin{proof}
Relative to the flag $F$,
$X_{\sigma[2]}\modmod S = \Pi_{\Hom(V,Q)}(Y_\sigma) \cap \Hom(V/S,Q)$ and
$X_{\sigma(01)} = Y_\sigma \cap \Hom(V/S,Q)$.  These are the same
as $\Pi_{\Hom(V,Q)}$ restricted to $\Hom(V/S,Q)$ is the identity map.
\end{proof}

\section{Proof of Horn's conjecture}
\label{sec:proof}

\subsection{The two-step flag manifold as a fibration}
Consider the map $q: Fl(d,r,\CC^n) \to Gr(r,\CC^n)$ defined by
$q(S,V) = V$.  This is a fibration, and the fibre over $V$ is
simply the Grassmannian $Gr(d,V)$.

Schubert varieties behave well under this fibration.  Let $F$ be
a flag on $\CC^n$, and $\sigma$ a $012$-string representing a 
Schubert class on $Fl(d,r, \CC)$.  
Then the restriction of $q$ to the Schubert cell $\Omega^F_\sigma$ 
is also a fibration.  We have
$q(\Omega^F_\sigma) = \Omega^F_{\sigma[2]}$ and the fibre over $V$
is the Schubert cell $\Omega^{F_V}_{\sigma(12)}$.  
Thus from this picture, we see that
intersection of Schubert varieties in 
$Fl(d,r, \CC^n)$ is related to two Schubert intersection problems in
Grassmannians: one on the base space $Gr(r,\CC^n)$ and one on
the fibre $Gr(d,V)$.  The precise relationship is as follows:

\begin{lemma}
Let $F^1, \ldots, F^s$ be generic flags on $\CC^n$.
Assume that $S^{\sigma^1(12)} \cdots S^{\sigma^s(12)} \neq 0 \in H^*(Gr(d,r))$.
Then $\bigcap_{i=1}^s \Omega^{F^i}_{\sigma^i}$ has a point of intersection
if and only if $\bigcap_{i=1}^s \Omega^{F^i}_{\sigma^i[2]}$ has 
a point of intersection.
\end{lemma}

\begin{proof}
The ``only if'' direction is clear: if 
$\bigcap_{i=1}^s \Omega^{F^i}_{\sigma^i}$ is non-empty, then so is 
its image under $q$.  

Thus suppose 
$I = \bigcap_{i=1}^s \Omega^{F^i}_{\sigma^i[2]}$ is non-empty.
$S^{\sigma^1(12)} \cdots S^{\sigma^s(12)} \neq 0 \in H^*(Gr(d,r))$ is equivalent to
the Schubert problem in the fibre
$\bigcap_{i=1}^s \Omega^{F_V^i}_{\sigma^i(12)}$ having a point of
intersection for $V$ which induce generic flags $F^i_V$.  But
by remark \ref{rmk:genericinducedflags}, a generic $V \in I$ has
this property.  Thus there is a point in
$\bigcap_{i=1}^s \Omega^{F^i}_{\sigma^i}$ over a generic point in
$I$.  In particular this intersection is non-empty.
\end{proof}

This simple lemma gives us a way (in fact many different ways) of
turning a Grassmannian Schubert intersection problem into a
Schubert intersection problem
on a two-step flag manifold.  Note that given any $01$-strings 
$\tau$ and $\mu$ with the correct number of \zeros and \ones, we can 
always construct a $012$-string $\sigma$ such that $\sigma[2]=\tau$ and 
$\sigma(12) = \rho$. 

\begin{definition}
We call $(\sigma^1, \ldots, \sigma^s)$ a {\bf lifting} of
$({\tau^1}, \ldots, {\tau^s})$, if $\sigma^i[2] = \tau^i$,
$\sigma^i(12) = \rho^i$ and
$S^{\rho^1} \cdots S^{\rho^s} \neq 0 \in H^*(Gr(d,r))$.
\end{definition}

Thus we see that there is a lifting corresponding to each 
$s$-tuple $(\rho_1, \ldots, \rho_s)$ such that
$S^{\rho^1} \cdots S^{\rho^s} \neq 0 \in H^*(Gr(d,r))$.
If we allow $d$ to vary between $1$ and $r$, we also have one Horn 
inequality for each such triple.  This is no
coincidence: we are about to see that each Horn inequality
arises from applying corollary \ref{cor:almosthorn} to a lifting.

\subsection{Necessity of the Horn inequalities}
We show that each lifting of a Grassmannian Schubert problem gives
rise to a Horn inequality.  This argument will establish the
necessity of the Horn inequalities. 

Suppose $\lambda^1, \ldots, \lambda^s \in \Lambda(r,n-r)$, with
$$S^{\lambda^1} \cdots S^{\lambda^s} \neq 0 \in H^*(Gr(r,\CC^n))$$
Then $\bigcap_{i=1}^s \Omega^{F^i}_{\lambda^i}$ contains a point
of intersection. Thus so does any lifting of this Schubert problem.

Let $(\sigma^1, \ldots, \sigma^s)$ be such a lifting, say, corresponding
to an $s$-tuple of $01$-strings $(\rho^1, \ldots, \rho^s)$.  
Let $\mu^i$ denote the partition corresponding to $\rho^i$.
Since
$\bigcap_{i=1}^s \Omega^{F^i}_{\sigma^i}$ is non-empty, for generic
flags, the intersection of tangent spaces
$\bigcap_{i=1}^s Y_\sigma^i$ must be transverse.  By 
corollary \ref{cor:almosthorn} this means that
\begin{equation}
\label{eqn:althorn}
\sum_{i=1}^s \dim X_{\sigma^i(02)} \geq (s-1)d(n-r)
\end{equation}
is a necessary inequality.

We now compute $\dim X_{\sigma(02)}$ (for ease of notation we are fixing
$i$ and omitting the superscript $i$ from $\sigma, \rho, \mu, \lambda$).
Let $z(j)$ denote the number
of \zeros in the list $\sigma_1, \cdots, \sigma_{j-1}$.
\begin{align*}
\dim X_{\sigma(02)} &= \#\{(j',j)\ |\ j'<j,\ \sigma_{j'}=0,\ \sigma_j=2\} \\
&= \sum_{j\ |\ \sigma_j =2} \lambda_{j - z(j)} \\
&= \sum_{k=1}^d \lambda_{\text{position of $k^{\rm th}$ \one in $\rho$}} \\
&= \sum_{k=1}^d \lambda_{\mu_k+k}
\end{align*}

Rewriting (\ref{eqn:althorn}) based on this calculation,
we obtain
$$
\sum_{i=1}^s \sum_{k=1}^d \lambda^i_{\mu^i_k +k} \geq (s-1)d(n-r).
$$
which is exactly the inequality (\ref{eqn:hornineq}).
Thus we see that the inequality (\ref{eqn:althorn}) is precisely the
Horn inequality corresponding to $(\mu_1, \ldots, \mu_s)$.

\subsection{Sufficiency of the Horn inequalities}
To prove sufficiency of the Horn inequalities, we must show that
whenever $S^{\lambda^1} \cdots \lambda^s = 0 \in H^*(Gr(r,\CC^n))$,
there is a Horn inequality violated by the $\lambda^i$.  We
have already seen that a violation of the Horn inequalities gives
rise to a non-transverse intersection 
$\bigcap_{i=1}^s X_{\sigma^i(02)}$ for
a lifting $(\sigma^1, \ldots, \sigma^s)$ of 
$(\lambda^1, \ldots, \lambda^s)$.  We would now like to prove the
reverse: a nontransverse intersection 
$\bigcap_{i=1}^s X_{\sigma^i(02)}$ leads to a violation of a Horn
inequality.  This requires an inductive argument.

\subsubsection{The inductive step}
\begin{lemma}
\label{lem:induction}
Let $\lambda^1, \ldots, \lambda^s \in \Lambda(r,n-r)$, and
let $(\sigma^1, \ldots, \sigma^s)$ be a lifting of 
$(\lambda^1, \ldots, \lambda^s)$ to Schubert varieties in 
$Fl(d,r,\CC^n)$.  Let $\mu^i$
be the partition corresponding to $\sigma^i(02)$.
If $\lambda^1, \ldots, \lambda^s$ satisfies all of its
the Horn inequalities, then $\mu^1, \ldots, \mu^s$ must
satisfy all of its Horn inequalities.
\end{lemma}

This lemma allows us to argue by induction.  Suppose the
Horn conditions are sufficient for all integers $d$, with $d<r$.
Suppose $S^{\lambda^1} \cdots S^{\lambda^s} = 0 \in H^*(Gr(r,\CC^n))$.
We show that one of the two things must be true.
\begin{enumerate}
\item The product is zero in cohomology for dimensional reasons.
\begin{center} \emph{or} \end{center}
\item There is a lifting of $(\lambda^1, \ldots, \lambda^s)$
to $(\sigma^1, \ldots, \sigma^s)$ such that the intersection
$\bigcap_{i=1}^s X_{\sigma^i(02)}$ is non-transverse.
\end{enumerate}
In the first case, the Horn inequality for $d=r$ is violated.
In the second case, our inductive hypothesis tells us that
some Horn inequality is violated by 
$\sigma^1(02), \ldots, \sigma^s(02)$.  So by 
lemma \ref{lem:induction}, there must be some Horn inequality violated 
by $\lambda^1, \ldots, \lambda^s$.

\begin{proof}
Let $V \subset \CC^n$, with $\dim V = r$ and
fix flags $F^1, \ldots, F^s$ on $\CC^n$ which are almost generic
for $V$.
Let $S \subset V$, with $\dim S = d$, and
$S$ in Schubert position $\sigma^i(12)$ with respect to
the induced flag $F^i_V$.

We work inside the two-step flag manifold $Fl(d',d,V)$.
let $(\chi^1, \ldots, \chi^s)$ be a lifting of 
$(\sigma^1(12), \ldots, \sigma^s(12))$.  Of course, since
this is a lifting 
$$S^{\chi^1(12)} \cdots S^{\chi^s(12)} \neq 0 \in H^*(Gr(d',S)).$$
Also, the intersection of
Schubert varieties
$\bigcap_{i=1}^s \Omega^{F^i_V}_{\sigma^i(12)}$ is non-empty, since
it contains the point $S$.  Thus 
$\bigcap_{i=1}^s \Omega^{F^i_V}_{\chi^i}$ is non-empty.

Now we consider the fibration $p:Fl(d',d,V) \to Gr(d',V)$.  As
in the case with $q$, the fibration $p$ maps 
Schubert varieties map to Schubert varieties.  The image is
$p(\Omega^{F}_\chi) = \Omega^F_{\chi[1]}$.  Thus we see that
$$S^{\chi^1[1]} \cdots S^{\chi^s[1])} \neq 0 \in H^*(Gr(d',V)).$$

We now check that the Horn inequality for $\mu^1, \ldots, \mu^s$
corresponding to $(\chi^1(12), \ldots, \chi^s(12))$ is identical
to the Horn inequality for $\lambda^1, \ldots, \lambda^s$, corresponding
to $(\chi^1[1], \ldots, \chi^s[1])$.  The two inequalities are
\begin{equation}
\label{eqn:lambdaeqn}
\sum_{i=1}^s \sum_{k=1}^d 
\lambda^i_{\text{position of the $k^{\rm th}$ \one in $\chi^i[1]$}}
\geq (s-1)d'(n-r)
\end{equation}
and
\begin{equation}
\label{eqn:mueqn}
\sum_{i=1}^s \sum_{k=1}^d 
\mu^i_{\text{position of the $k^{\rm th}$ \one in $\chi^i(12)$}}
\geq (s-1)d'(n-r).
\end{equation}
Now, $\lambda^i_l = \mu^i_{l'}$ where $l-l'$ is the number of
\zeros in the list $\chi^i_1, \cdots, \chi^i_{l-1}$, thus $l'$.
But $\chi^i[1]_l = 1$ if and only if $\chi^i_l = 2$ if and only if
$\chi^i(12)_{l'} = 1$.  Thus we see that
$$
\lambda^i_{\text{position of the $k^{\rm th}$ \one in $\chi^i[1]$}}
=\mu^i_{\text{position of the $k^{\rm th}$ \one in $\chi^i(12)$}}
$$
and the two Horn inequalities (\ref{eqn:lambdaeqn}) and 
(\ref{eqn:mueqn}) are the same.

Since every Horn inequality for
$\mu^1, \ldots, \mu^s$ arises in this way, if 
$(\lambda^1, \ldots, \lambda^s)$ satisfies all its Horn inequalities,
then so does $(\mu^1, \ldots, \mu^s)$.
\end{proof}

\begin{remark} It is perhaps most natural to view this lemma
as a statement about the three-step flag manifold 
$Fl(d',d,r,\CC^n)$.  We take a Schubert problem on $Gr(r,\CC^n)$
and lift it to $Fl(d,r,\CC^n)$.  We then lift this again, to
a problem on $Fl(d',d,r,\CC^n)$, given by $0123$-strings
$\omega^1, \ldots, \omega^s$.  We can then interpret both
inequalities (\ref{eqn:lambdaeqn}) and (\ref{eqn:mueqn}), as the
statement $\sum_{i=1}^s \dim X_{\sigma^i(03)} \geq (s-1)d(n-r)$,
which can be viewed as a Horn inequality for both $\lambda^i$
and $\mu^i$.
\end{remark}

\subsubsection{Proof of sufficiency}
We now have all the ingredients in place to prove the sufficiency
of the Horn conditions.

\begin{proof}(Horn's conjecture)
As always, let $V \subset \CC^n$, with $\dim V = r$, and fix
almost generic flags $F^i$.
Consider a generic $\phi \in \bigcap_{i=1}^s X_\lambda^i$.  Let
$S$ be the kernel of $\phi$ in Schubert position $\rho^i$ with
respect to generic flags $F^i_V$ on $V$.  Let 
$(\sigma^1, \ldots, \sigma^s)$ be the lifting of 
$(\lambda^1, \ldots, \lambda^s)$ by $(\rho^1, \ldots, \rho^s)$.

Assume that
$S^{\lambda^1} \cdots S^{\lambda^s} = 0 \in H^*(Gr(r,\CC^n))$.
Thus
$\bigcap_{i=1}^s Y_{\sigma^i}$
is non-transverse.

There are two possibilities.  Either $S=V$ or $\dim S < \dim V$.
In the first case, we must have $\phi=0$, which means that
$\bigcap_{i=1}^s X_\lambda^i = \zeroset$ (otherwise, $\phi=0$ would
not be a generic choice).  If 
$\sum_{i=1}^s \codim X_\lambda^i = r(n-r)$ then this is a transverse
intersection, contradicting $S^{\lambda^1} \cdots S^{\lambda^s} = 0$.
Thus $\sum_{i=1}^s \codim X_\lambda^i > r(n-r)$ which violates
the Horn inequality for $d=r$.

If $\dim S< \dim V$,
we note the following facts:
$$\bigcap_{i=1}^s X_{\sigma^i(12)} =
\bigcap_{i=1}^s X_{\rho^i}$$
is a transverse intersection, since the Schubert varieties
$\Omega^{F^i_V}_{\rho_i}$ are in general position, and contain
$S$ as a point of intersection.  Also
$$\bigcap_{i=1}^s X_{\sigma^i(01)} =
\bigcap_{i=1}^s X_{\sigma^i[1]}\modmod S
= \bigcap_{i=1}^s X_{\lambda^i}\modmod S$$
is a transverse intersection, by corollary 
\ref{cor:transversereduction}.

If $\bigcap_{i=1}^s X_{\sigma^i(02)}$ is also transverse, then
by lemma \ref{lem:splitting}, $\bigcap_{i=1}^s Y_{\sigma^i}$ would
be transverse.  Thus it is not a transverse intersection, and
so by the inductive argument, some Horn inequality is violated.
\end{proof}

\subsection{Examples}

The proof of sufficiency describes a procedure for finding a Horn
inequality which is violated, when
$S^{\lambda^1} \cdots S^{\lambda^s} = 0.$
We now give some examples to show what happens in this process,
in the simplest case which is $s=2$.  (The sufficiency of the
Horn inequalities is an easy fact when $s=2$; nevertheless, it
illustrates the method of the proof fairly 
adequately.)

\begin{example}
\label{ex:simple}
Let 
$$\lambda^1
= {\bf 0} \leq 0 \leq 3 \leq 3 \leq {\bf 4}$$
and 
$$\lambda^2
= {\bf 0} \leq 1 \leq 3 \leq 3 \leq {\bf 4}$$
As in example \ref{ex:oppositeintersection}, we can illustrate
$X_{\lambda^1}$ and $X_{\lambda^2}$ in general position as
$$
X_{\lambda^1} = 
\left (
\begin{matrix}
0 & * & *  \\
0 & * & *  \\
0 & * & *  \\
0 & 0 & 0  
\end{matrix}
\right )
\qquad
X_{\lambda^2} = 
\left (
\begin{matrix}
 0 & 0 & 0 \\
 + & + & 0 \\
 + & + & 0 \\
 + & + & + 
\end{matrix}
\right )
$$
We'll write these both in a single diagram as
$$
\begin{tabular}{ccc} \hline
\sqempty   &\squp   & \squpend  \\\hline
\sqdown & \squpdown & \squpend  \\\hline
\sqdown & \squpdown & \squpend  \\\hline
\sqdown & \sqdown & \sqdownend \\\hline
\end{tabular}
$$
Take a generic point $\phi \in X_{\lambda^1} \cap X_{\lambda^2}$,
e.g.
$$ \phi = 
\left (
\begin{matrix}
0 & 0 & 0 \\
0 & 3 & 0 \\
0 & 8 & 0 \\
0 & 0 & 0 \\
\end{matrix}
\right )
$$
The kernel of $\phi$ has Schubert position $(101, 101)$.  We
use this Schubert position to lift $(\lambda^1, \lambda^2)$:
$$
\begin{tabular}{ccc} \hline
 \squp   & \sqempty   & \squpend \\\hline
 \squpdown &  \sqdown  & \squpend \\\hline
 \squpdown &  \sqdown  & \squpend \\\hline
 \sqdown   &  \sqdown  & \sqdownend \\\hline 
      & \sqdown & \squpend \\\cline{2-3}
\end{tabular}
$$
The upper right block is non-transverse for dimensional reasons.
Thus we are led to consider the Horn inequality corresponding to 
$(101,101)$,
i.e. 
$$\lambda^1_1 + \lambda^1_3 + \lambda^2_1 + \lambda^2_3 \geq 8.$$
The positions of the \ones in $(101,101)$ are the indices which
appear to the indices which appear on the left hand side.  We
see that this inequality is violated by $(\lambda^1, \lambda^2)$,
as 
$$\lambda^1_1 + \lambda^1_3 + \lambda^2_1 + \lambda^2_3 \geq 
0+3+1+3 <8.$$
\end{example}

\begin{example}
Let $$\lambda^1 
= {\bf 0} \leq 0 \leq 2 \leq 3 \leq 3 \leq 3 \leq 4 \leq {\bf 4}$$
and 
$$\lambda^2
= {\bf 0} \leq 1 \leq 1 \leq 3 \leq 3 \leq 3 \leq 3 \leq {\bf 4}.$$
We illustrate $X_{\lambda^1}$ and $X_{\lambda^2}$ as:
$$
\begin{tabular}{cccccc} \hline
\sqempty   & \squp   & \squp   & \squp   & \squp & \squpend \\\hline
\sqdown & \squpdown & \squpdown & \squpdown & \squp & \squpend \\\hline
\sqdown & \sqdown   & \squpdown & \squpdown & \squp & \squpend \\\hline
\sqdown & \sqdown   & \sqdown   & \sqdown   & \sqdown & \squpdownend \\\hline
\end{tabular}
$$
Let $\phi \in X_{\lambda^1} \cap X_{\lambda^2}$ be a generic
element, e.g.
$$ \phi = 
\left (
\begin{matrix}
0 & 0 & 0 & 0 & 0 & 0 \\
0 & 5 & 6 & 7 & 0 & 0 \\
0 & 0 & 8 & 9 & 0 & 0 \\
0 & 0 & 0 & 0 & 0 & 1 
\end{matrix}
\right )
$$
The kernel of $\phi$ has Schubert position $(100110,010011)$.  We
use this Schubert position to lift $(\lambda^1, \lambda^2)$:
$$
\begin{tabular}{cccccc} \hline
  \squp   & \squp   & \squp  & \sqempty   &\squp   & \squpend \\\hline
  \squpdown & \squpdown & \squp  & \sqdown &\squpdown & \squpend \\\hline
  \sqdown   & \squpdown & \squp  & \sqdown &\squpdown & \squpend \\\hline
  \sqdown   & \sqdown   & \squpdown& \sqdown &\squpdown & \sqdownend  \\\hline
      & &  & \sqdown & \squpdown  & \squpend \\\cline{4-6}
 &  &  & \sqdown & \squpdown  & \squpend \\\cline{4-6}
 &  &  & \sqdown & \sqdown  & \sqdownend \\\cline{4-6}
\end{tabular}
$$
The upper left, and lower right blocks have a transverse intersection.
However, the upper right block 
$$
\begin{tabular}{ccc} \hline
\sqempty   &\squp   & \squpend  \\\hline
\sqdown & \squpdown & \squpend  \\\hline
\sqdown & \squpdown & \squpend  \\\hline
\sqdown & \sqdown & \sqdownend \\\hline
\end{tabular}
$$
does not.  This block represents $X_{{\lambda^1}'} \cap X_{{\lambda^2}'}$
for 
$${\lambda^1}'
= {\bf 0} \leq 0 \leq 3 \leq 3 \leq {\bf 4}$$
$${\lambda^2}'
= {\bf 0} \leq 1 \leq 3 \leq 3 \leq {\bf 4}$$
In example \ref{ex:simple} we found that the Horn inequality
${\lambda^1}'_1 + {\lambda^1}'_3 + {\lambda^2}'_1 + {\lambda^2}'_3 \geq 8$,
which corresponds to $(101,101)$ is violated.
To find
a corresponding Horn inequality which is violated by 
$(\lambda^1, \lambda^2)$ we lift the Schubert position of
$\ker \phi$ by $(101, 101)$ to get $(200120, 020012)$.  The positions
of the \twos give the indices which appear in the relevant inequality.
In this case we find that the Horn inequality
$$\lambda^1_1 + \lambda^1_5 + \lambda^2_2 + \lambda^2_6 \geq 8$$
is violated.  Indeed
$$\lambda^1_1 + \lambda^1_5 + \lambda^2_2 + \lambda^2_6  = 
0+3+1+3 <8.$$
\end{example}

\end{document}